\documentclass{ws-rv975x65}   
\usepackage{ws-rv-thm}    
\usepackage[square]{ws-rv-van}    
\usepackage{ws-crop975x65}     
\usepackage{bbold, mathrsfs}

\usepackage{ws-minitoc}        

\allowdisplaybreaks[4]

\numberwithin{equation}{section} 
\numberwithin{theorem}{section} 
\numberwithin{remark}{section} 
\numberwithin{proposition}{section} 
\numberwithin{lemma}{section} 
\numberwithin{corollary}{section} 
\numberwithin{definition}{section} 
\numberwithin{example}{section} 

\makeatletter
\renewcommand*\l@chapter[2]{}
\renewcommand*\l@author[2]{}
\makeatother

\begin{document}

\setcounter{chapter}{5}
\setcounter{tocdepth}{2}
\dominitoc

\chapter{Stein's method, Malliavin calculus, Dirichlet forms and the
fourth moment theorem}
\chaptermark{Stein's method and Dirichlet forms}

\setcounter{page}{107}  

\author[Louis H.Y. Chen and Guillaume Poly]
{Louis H.Y. Chen\footnote{Department of Mathematics, National University of Singapore,
10 Lower Kent Ridge Road, Singapore 119076, Republic of Singapore, \
matchyl@nus.edu.sg} and
Guillaume Poly\footnote{Mathematical research unit, University of Luxembourg,
6 rue Richard Coudenhove-Kalergi, L-1359 Luxembourg,
Grand Duchy of Luxembourg, \
poly@crans.org}}


\begin{abstract}
The fourth moment theorem provides error bounds of the order
$\sqrt{{\mathbb E}(F^4) - 3}$ in the central limit theorem for elements $F$ of
Wiener chaos of any order such that ${\mathbb E}(F^2) = 1$.
It was proved by Nourdin and Peccati \cite{NP09} using Stein's method and the
Malliavin calculus. It was also proved by Azmoodeh, Campese and Poly~\cite{A13}
using Stein's method and Dirichlet forms.  This paper is an exposition on the
connections between Stein's method and the Malliavin calculus and between Stein's
method and Dirichlet forms, and on how these connections are exploited in proving
the fourth moment theorem.
\end{abstract}

\smalltoc
\tableofcontents

\pagebreak  

\body
\vspace*{-35pt} 
\section{Introduction}

In 2005, Nualart and Peccati \cite{NP05}, discovered the remarkable fact that a sequence of multiple Wiener-It\^{o} integrals, that is, members of a Wiener chaos, converges in distribution to a Gaussian random variable if and only if their second and fourth moments converge to the corresponding moments of the limiting random variable. The proof in \cite{NP05} relies on a theorem in stochastic calculus, namely the Dubins-Schwarz Theorem. Although the proof is elegant, it does not provide good bounds on the distances (for instance, the Kolmogorov distance) between the sequence and its limit.

A few years later, the result of Nualart and Peccati \cite{NP05} was generalized and developed further by Peccati and Tudor \cite{PT04} and also by Nourdin and Peccati \cite{NP09}. Peccati and Tudor \cite{PT04} considered the multivariate central limit theorem for vectors of multiple Wiener-It\^{o} integrals and established that if the covariance matrices of the vectors of multiple Wiener-It\^{o} integrals converge to that of a Gaussian vector, then joint convergence in distribution to Gaussian is equivalent to coordinate-wise convergence in distribution to Gaussian. As an application of this result, the central limit theorem for any sequence of random variables can be established by proving a central limit theorem for each term in the chaos expansion.

Another significant development from the work of Nualart and Peccati \cite{NP05} is that of Nourdin and Peccati \cite{NP09}. Developing considerably a remarkable intuition of Nualart and Ortiz-Latorre in \cite{Nu08}, Nourdin and Peccati \cite{NP09} established a fundamental connection between Stein's method and the Malliavin calculus via the use of the Malliavin operators ($D, \delta, {\textbf L}$). This connection resulted in the {derivation} of errors bounds, often optimal, in the central limit theorems for random variables in the Wiener space. Of special interest in this paper is the {proof} of error bounds in the central limit theorem for multiple Wiener-It\^{o} integrals, which we call the fourth moment theorem.

The success of such a connection relies on the fact that both Stein's method and the Malliavin calculus are built on some integration by parts techniques. In addition, the operators of the Malliavin calculus, $D, \delta, {\textbf L}$, satisfy several nice integration by parts formulae which fit in perfectly with the so-called Stein equation. For a good overview of these techniques, we refer to the following website.

\begin{verbatim}https://sites.google.com/site/malliavinstein/home\end{verbatim}

The work of Nourdin and Peccati \cite{NP09} has added a new dimension to Stein's method. Their approach of combining Stein's method with the Malliavin calculus has led to improvements and refinements of many results  in probability theory, such as the Breuer-Major theorem \cite{BM83}. More recently, this approach has been successfully used to obtain central limit theorems in stochastic geometry, stochastic calculus, statistical physics, and for zeros of random polynomials. It has also been

\pagebreak  

\noindent   
extended to different settings as in non-commutative probability and Poisson chaos. {Of particular interest is the connection between the Nourdin-Peccati analysis and information theory, which was recently revealed in \cite{Le14,NOS13}.} An overview of these new developments can also be found in the above website.

The approach of Nourdin and Peccati \cite{NP09} entails the use of the so-called product formula for Wiener integrals. The use of this formula makes the proofs rather involved since it relies on subtle combinatorial arguments. Very recently, {starting with the work of Ledoux \cite{le}}, a new approach to the fourth moment theorem was developed by Azmoodeh, Campese and Poly \cite{A13} by combining Stein's method with the Dirichlet form calculus. An advantage of this new approach is that it provides a simpler proof of the theorem by avoiding completely the use of the product formula. {Moreover, since a Dirichlet space is a more general concept than the Wiener space, the former contains examples of fourth moment theorems that cannot be realized on the latter.}

 A more algebraic flavor of this approach has enabled Azmoodeh, Malicet and Poly \cite{A14} to prove that convergence of pairs of moments other than the 2nd and 4th (for example, the 6th and 68th) also implies the central limit theorem. This new approach seems to open up new possibilities and perhaps also central limit theorems on manifolds.

This paper is an exposition on the connections between Stein's method and the Malliavin calculus and between Stein's method and Dirichlet forms, and on how these connections are exploited in proving the fourth moment theorem.

\section{Stein's method}
\subsection{How it began}

Stein's method began with Charles Stein using his own approach in the 1960's to prove the combinatorial central limit theorems of Wald and Wolfowitz \cite{W44} and of Hoeffding \cite{H51}.  Motivated by permutation tests in nonparametric statistics, Wald and Wolfowitz \cite{W44} proved that under certain conditions, $\sum_{i=1}^n a_ib_{\pi(i)}$ converges in distribution to the standard normal distribution $\mathcal{N}(0,1)$, where $\{a_i, b_j: i, j = 1, \dots, n\}$ are real numbers and $\pi$ a random permutation of $\{1, \dots, n\}$. Hoeffding \cite{H51} generalized the result of Wald and Wolfowitz \cite{W44} to $\sum_{i=1}^n c_{i\pi(i)}$, where $\{c_{ij}: i,j = 1,\dots,n\}$ is a square array of real numbers and $\pi$ a random permutation of $\{1, \dots, n\}$.
\medskip

Let $W = \sum_{i=1}^n c_{i\pi(i)}$ and let $\phi$ be the characteristic function of $W$. Assume $c_{i\cdot}=c_{\cdot j}=0$ where $c_{i\cdot}=\sum_{j=1}^n c_{ij}/n$ and $c_{\cdot j}=\sum_{i=1}^n c_{ij}/n$ (which  {implies} ${\mathbb E}[W]=0$), and also assume ${\rm Var}(W) = 1$. While using exchangeable pairs to show that

\begin{equation}\label{Stein-complex}
\phi'(t)\simeq -t \phi(t)
\end{equation}
Stein realized that there was nothing special about the complex exponentials. Since $\phi'(t) = i{\mathbb E} \big[We^{itW}\big]$ and $-t\phi(t) = - t{\mathbb E} \big[e^{itW}\big] = - \frac{1}{i}{\mathbb E} \left[\frac{\partial}{\partial W}e^{itW}\right]$,
replacing the complex exponential by an arbitrary function $f$, (\ref{Stein-complex}) becomes
$$
{\mathbb E} \big[Wf(W)\big]\simeq {\mathbb E} \big[f'(W)\big].
$$
{By letting $f$ be a solution}, say $f_h$, of the differential equation
$$
f'(w) - wf(w) = h(w) - {\mathbb E} \big[h(Z)\big],
$$
where $h$ is a bounded function and $Z \sim \mathcal{N}(0,1)$, one obtains
$$
{\mathbb E} \big[h(W)\big] - {\mathbb E}\big[ h(Z)\big] = {\mathbb E}\big[f_h'(W) - Wf_h(W)\big].
$$
How close ${\mathscr L}(W)$ is to $\mathcal{N}(0,1)$ can then be determined by bounding ${\mathbb E}\big[f_h'(W) - Wf_h(W)\big]$. There is no inversion of the characteristic function.
\medskip

This story of how Stein's method began is based on personal communications with Charles Stein and Persi Diaconis and also on an interview
of Stein in Leong \cite{L10}.

\subsection{A general framework}
Stein's method for normal approximation was {published} in his seminal 1972 paper in the Proceedings of the Sixth Berkeley Symposium. Although the method was developed for normal approximation, Stein's
ideas were very general and the method was modified by Chen \cite{C75} for
Poisson approximation. Since then the method has been constantly developed and applied to many
approximations beyond normal and Poisson and in finite as well as infinite
dimensional spaces. As Stein's method works well for dependent random variables, it has been applied, and continues to be applied, to a large number of problems in many different fields. The method, together with its applications, continues to grow and remains a very active research area.
See, for example, Stein \cite{S86}, Arratia, Goldstein and Gordon \cite{A90}, Barbour, Holst and Janson \cite{BHJ92}, Diaconis and Holmes \cite{D04}, Barbour and Chen \cite{BC05a,BC05b}, Chatterjee, Diaconis and Meckes \cite{C05}, Chen, Goldstein and Shao \cite{C11}, Ross \cite{R11}, Shih \cite{S11}, Nourdin and Peccati \cite{N12},
and Chen and R\"ollin \cite{C13a,C13b}.

In a nutshell, Stein's method can be described as follows. Let $W$ and $Z$ be random elements taking values in a space $\mathcal{S}$ and let $\mathcal{X}$ and $\mathcal{Y}$ be
some classes of real-valued functions defined on $\mathcal{S}$. In approximating the distribution ${\mathscr L}(W)$ of $W$ by the distribution ${\mathscr L}(Z)$ of $Z$, we
write ${\mathbb E} \big[h(W)\big] -  {\mathbb E} \big[h(Z)\big] = {\mathbb E} \big[Lf_h(W)\big]$ for a test function $h\in \mathcal{Y}$,
where $L$ is a linear operator (Stein operator) from $\mathcal{X}$ into $\mathcal{Y}$ and $f_h \in \mathcal{X}$ a
solution of the equation
\begin{equation}
  Lf = h - {\mathbb E} \big[h(Z)\big]\qquad {\text{(Stein equation)}.}
\end{equation}

The error
${\mathbb E}\big[Lf_h(W)\big]$ can then be bounded by studying the solution $f_h$ and
exploiting the probabilistic properties of $W$. The operator $L$ characterizes ${\mathscr L}(Z)$  in the sense that ${\mathscr L}(W) = {\mathscr L}(Z)$ if and only if for a
sufficiently large class of functions $f$ we have
\begin{equation}
  {\mathbb E}\big[Lf(W)\big]= 0 \qquad {\text{(Stein identity)}.}
\end{equation}

In normal approximation, where ${\mathscr L}(Z)$ is the standard normal distribution,
the operator used by Stein \cite{S72} is given by $Lf(w) = f'(w) - wf(w)$ for
$w\in \mathbb{R}$, and in Poisson approximation, where ${\mathscr L}(Z)$ is the
Poisson distribution with mean $\lambda > 0$, the operator $L$ used by Chen \cite{C75}
is given by $Lf(w) = \lambda f(w + 1) - wf(w)$ for $w \in \mathbb{Z}_+$.
However the operator $L$ is not unique even for the same approximating distribution
but depends on the problem at hand.
For example, for normal approximation $L$ can also be taken to be the generator of
the Ornstein-Uhlenbeck process, that is, $Lf(w) = f''(w) - wf'(w)$, and for Poisson
approximation, $L$ taken to be the generator of an immigration-death process,
that is, $Lf(w) = \lambda[f(w + 1) - f(w)] + w[f(w - 1) - f(w)]$.
This generator approach, which is due to Barbour \cite{B88}, allows
extensions to multivariate and process settings. Indeed, for multivariate normal approximation, $Lf(w) = \Delta f(w) - w\cdot
\nabla f(w)$, where $f$ is defined on the Euclidean space; see
Barbour \cite{B90} and G\"otze \cite{G91}, and also Reinert and R\"ollin \cite{RR09} and Meckes \cite{M09}.

\subsection{Normal approximation}

In many problems of normal approximation, the random variable $W$ whose distribution is to be approximated satisfies this equation
\begin{equation}\label{Stein-general}
{\mathbb E} \big[Wf(W)\big] = {\mathbb E} \big[T_1 f'(W + T_2)\big].
\end{equation}
where $T_1$ and $T_2$ are some random variables defined on the same probability space as $W$, and $f$ is an absolutely continuous function for which the expectations in (\ref{Stein-general}) exist. Examples of $W$ satisfying this equation include sums of locally dependent random variables as considered in Chen and Shao \cite{CS04} and exchangeable pairs as defined in Stein \cite{S86}. More generally, a random variable $W$ satisfies (\ref{Stein-general}) if there is a Stein coupling $(W,W',G)$ where $W, W',G$ are defined on a common probability space such that ${\mathbb E} \big[Wf(W)\big]= {\mathbb E}\big[Gf(W') - Gf(W)\big]$ for absolutely continuous functions $f$ for which the expectations exist (see {Chen and R\"ollin \cite{C13b}}). In all cases it is assumed that ${\mathbb E}\big[W\big]= 0$ and ${\rm Var} (W) = 1$. Letting $f(w) = w$, we have $1 = {\mathbb E} \big[W^2\big] = {\mathbb E} \big[T_1\big].$

As an illustration, let $(W,W')$ be an exchangeable pair of random variables with ${\mathbb E} (W) = 0$ and $\mathrm{Var}(W) = 1$ such that ${\mathbb E} \big[W' - W|W\big] = -\lambda W$ for some $\lambda > 0$. Since the function $(w,w') \longmapsto (w'-w)(f(w') + f(w))$ is anti-symmetric, the exchangeability of $(W,W')$ implies
$$
E\big[(W'-W)(f(W')+f(W))\big] = 0.
$$
From this we obtain
$$
{\mathbb E} \big[Wf(W)\big]=\frac{1}{2\lambda} {\mathbb E}\Big[(W'-W)(f(W')-f(W))\Big]
$$
$$
= \frac{1}{2\lambda}{\mathbb E}\Big[(W'-W)^2\int_0^1f'(W+ (W'-W)t)dt\Big] = {{\mathbb E} \big[T_1f'(W+T_2)\big]}
$$
where ${T_1 = {\displaystyle \frac{1}{2\lambda}(W'-W)^2}}$, ~$T_2 = (W'-W)U$, and $U$ uniformly distributed on $[0,1]$ and independent of ${W, W'}, ~T_1$ and $T_2$.

Let $f_h$ be the unique bounded solution of the Stein equation
\begin{equation}\label{Stein-eqn-normal}
 f'(w) - wf(w) = h(w) - {\mathbb E} \big[h(Z)\big],
\end{equation}

where $w \in {\mathbb R}$, ~$Z \sim \mathcal{N}(0,1)$, and $h$ a bounded test function.

The following boundedness properties of $f_h$ are useful for bounding the errors in the approximation.

If $h$ is bounded, then
$$
\|f_h\|_{\infty} \le \sqrt {2\pi}\|h\|_{\infty}, ~~\|f'_h\|_{\infty} \le
4\|h\|_{\infty}.
$$

If $h$ is absolutely continuous, then
$$
\|f_h\|_{\infty} \le 2\|h'\|_{\infty},~~\|f'_h\|_{\infty} \le \sqrt {2/\pi}\|h'\|_{\infty}, ~~\|f''_h\|_{\infty} \le 2\|h'\|_{\infty}.
$$

If $h = { \textbf{1}_{(-\infty,x]}}$, then for all $w, v \in {\mathbb R}$,
$$
0 \le f_h(w) \le \sqrt{2\pi}/4, ~~|wf_h(w)| \le 1, ~~|f_h'(w)| \le 1, ~~|f_h'(w) - f_h'(v)| \le 1.
$$

These can be found in Lemmas 2.4 and 2.5 of Chen, Goldstein and Shao \cite{C11} and in Lemma 2.2 of Chen and Shao \cite{CS05}

Assume that ${\mathbb E} \big[W\big] = 0$ and ${\rm Var} (W) = 1$. From (\ref{Stein-general}) and (\ref{Stein-eqn-normal}),
\begin{eqnarray}\notag
{\mathbb E} \big[h(W)\big] - {\mathbb E} \big[h(Z)\big] &=& {\mathbb E} \big[f'_h(W) - T_1f'_h(W + T_2)\big]\\ \notag
&=&{\mathbb E} \big[T_1(f'_h(W) - f'_h(W + T_2))\big] + {\mathbb E} \big[(1 - T_1)f_h'(W)\big]\\
\label{Stein-solu-property}
\end{eqnarray}

Different techniques have been developed for bounding the error term on the right side of (\ref{Stein-solu-property}). We will consider two special cases.

Case 1. Assume that $T_1 = 1$. This is the case of zero-bias coupling. See Goldstein and Reinert \cite{GR97}, and also Chen, Goldstein and Shao \cite{C11}.
Let $h$ be absolutely continuous such that $h'$ is bounded.  From (\ref{Stein-solu-property}),
$$
|{\mathbb E} \big[h(W)\big] - {\mathbb E} \big[h(Z)\big]| = {|{\mathbb E}\Big[\int_0^{T_2} f''_h(W +t )dt\Big]}| \le \|f''_h\|_{\infty} {\mathbb E} \big[|T_2|\big] \le 2\|h'\|_{\infty}{\mathbb E}\big[|T_2|\big].
$$
From this we obtain the following bound on the Wasserstein distance between ${\mathscr L}(W)$ and $\mathcal{N}(0,1)$.
$$
d_{\mathrm{W}}({\mathscr L}(W),\mathcal{N}(0,1)): = \sup_{|h(x) - h(y)| \le |x - y |} |{\mathbb E}\big[ h(W)\big] - {\mathbb E}\big[ h(Z)\big]| \le 2{\mathbb E}\big[|T_2|\big].
$$
Note that $d_{\mathrm{W}}({\mathscr L}(W),\mathcal{N}(0,1)) = \|F - \Phi\|_1$ where $F(x) = P(W \le x)$ and $\Phi(x) = P(Z \le x)$.

If $W = X_1 + \dots + X_n$ where $X_1 , \dots, X_n$ are independent random variables with ${\mathbb E} [X_i] = 0$, ~${\rm Var}(X_i)= \sigma^2$ and ${\mathbb E} \big[|X_i|^3\big] = \gamma_i < \infty$, then $T_2 = \xi_I - X_I$ where the $X_i$, the $\xi_i$ and $I$ are independent, ${\mathbb E} \big[|\xi_i|\big] = \gamma_i/2\sigma_i^2$, and $P(I=i) = \sigma_i^2$. Note that $\sum \sigma_i^2 = {\rm Var}(W) = 1$. Then the bound on the Wasserstein distance between ${\mathscr L}(W)$ and {$\mathcal{N}(0,1)$} is $2{\mathbb E}\big[|\xi_I - X_I|\big] \le 3\sum \gamma_i$.

It is more difficult to obtain a bound on the Kolmogorov distance between ${\mathscr L} (W)$
and ${\mathcal{N}(0,1)}$, namely $\sup_{x\in{\mathbb R}} |P(W \le x) - \Phi(x)|$.
Such a bound can be obtained by induction or the use of a concentration inequality.
For induction, see \cite{B84}. For the use of a concentration inequality,
see Chen \cite{C98} and Chen and Shao \cite{CS01}
for sums of independent random variables, and Chen and Shao \cite{CS04}
for sums of locally dependent random variables.
See also Chen, Goldstein and Shao \cite{C11}.
For sums of independent random variables, Chen and Shao \cite{CS01} obtained a bound of
$4.1\sum \gamma_i$ on the Kolmogorov distance.

Case 2. Assume that $T_2 = 0$. This is the case if $W$ is a functional of independent Gaussian random variables as considered by Chatterjee \cite{C09} or a functional of Gaussian random fields as considered by Nourdin and Peccati \cite{NP09}. In this case, (\ref{Stein-solu-property}) becomes
$$
{\mathbb E} \big[h(W)\big] - {\mathbb E} \big[h(Z)\big] = {\mathbb E} \big[(1 - T_1)f_h'(W)\big] = {\mathbb E}\big[(1 - {\mathbb E}[T_1|W])f_h'(W)\big].
$$
Let $h$ be such that $|h| \le 1$. Then we obtain the following bound on the  total variation distance between ${\mathscr L}(W)$ and $\mathcal{N}(0,1)$.

\begin{eqnarray*}
d_{\mathrm{TV}}({\mathscr L}(W), \mathcal{N}(0,1)) &:=& \frac{1}{2}\sup_{|h|\le 1}|{\mathbb E} \big[h(W)\big] - {\mathbb E} \big[h(Z)\big]| \\
&\le& \frac{1}{2}\|f_h'\|_{\infty}{\mathbb E}\big[|1 - {\mathbb E}[T_1|W]|\big]\\
&\le& 2\sqrt{{\rm Var}({\mathbb E}[T_1|W])}
\end{eqnarray*}

While Chatterjee \cite{C09} used a Poincar\'e inequality of second order to bound $2\sqrt{{\rm Var}({\mathbb E}[T_1|W])}$, Nourdin and Peccati \cite{NP09} deployed the Malliavin calculus. In the next two sections, we will discuss how the Malliavin calculus is used to bound $2\sqrt{{\rm Var}({\mathbb E}[T_1|W])}$.

\section{Malliavin calculus}

\subsection{A brief history}

The Malliavin calculus was born in 1976 in a symposium in Kyoto. Paul Malliavin presented a remarkable theory which extended the powerful calculus of variations to the framework of stochastic calculus. The initial goal of the theory was to provide a probabilistic proof of the H\"{o}rmander criterion (H\"{o}rmander \cite{H67})  of hypoellipticity by relating the smoothness of the solutions of a second order partial differential equation with the smoothness of the law of the solution of a stochastic differential equation. In order to prove that a random variable $X$ has a smooth density, Paul Malliavin introduced the following very efficient criterion.
\begin{lemma}
Assume that for each $k\ge 1$, there is $C_k>0$ such that for any $\phi\in\mathcal{C}^\infty_C$, the class of $\mathcal{C}^\infty$ functions with compact support, we have
\begin{equation}\label{Malliavinlma}
\Big|{\mathbb E}\left[\phi^{(k)}(X)\right]\Big|\le C_k \|\phi\|_\infty,
\end{equation}
then the distribution of $X$ has a $\mathcal{C}^\infty$ density.
\end{lemma}
In order to prove the inequality (\ref{Malliavinlma}) for $X$, Malliavin showed, by integration by parts, that for some suitable weight $H_k$,
\begin{equation}\label{IPP-Malliavin}
{\mathbb E}\big[\phi^{(k)}(X)\big]={\mathbb E}\big[\phi(X)H_k\big].
\end{equation}
Hence, (\ref{Malliavinlma}) holds with $C_k={\mathbb E}[|H_k|]$. In many situations of interest, $X=F(G_1,\cdots,G_n)$ where the $\{G_i\}_{i\ge 1}$ are i.i.d. Gaussian and $F$ is smooth. The equation (3.2) can be proved by integration by parts through the "Gaussian structure" of $X$ together with some "non-degeneracy" assumption on $F$. Since its introduction, the Malliavin calculus has been extended and used in many different areas of probability theory. However, regardless of the application, the central role of the Malliavin calculus always consists {in} proving that some integration by parts formula holds. In the present exposition, we use integration by parts not for proving the smoothness of a density, but for establishing Stein's bounds.

\subsection{Malliavin derivatives}

Let ${\mathscr H}$  be a real separable Hilbert space, typically ${\mathscr H}=L^2({\mathbb R}_+)$ but the particular choice of ${\mathscr H}$ does not matter. We denote by $X=\{X(h), h\in {\mathscr H}\}$ an \textit{isonormal Gaussian process} over ${\mathscr H}$. That means, $X$ is a centered Gaussian family of random variables defined in some probability space $(\Omega, \mathcal{F},P)$, with covariance given by
\[
{\mathbb E}[X(h)X(g)]= \langle h,g \rangle_{{\mathscr H}},
\]
for any $h,g \in {\mathscr H}$. We also assume that $\mathcal{F}$ is generated by $X$.
~\\~\\
Let $\mathcal{S}$ be the set of all cylindrical random variables of the form:
\begin{equation}
F=g\left( X(\phi _{1}),\ldots ,X(\phi _{n})\right) ,  \label{v3}
\end{equation}
where $n\geq 1$, $g:\mathbb{R}^{n}\rightarrow \mathbb{R}$ is an infinitely
differentiable function such that its partial derivatives have polynomial growth, and $\phi _{i}\in {\mathscr H}$,
$i=1,\ldots,n$.
The {\it Malliavin derivative}  of $F$ with respect to $X$ is the element of
$L^2(\Omega ,{\mathscr H})$ defined as

\begin{equation}\label{defi-malliavin}
DF\;=\;\sum_{i=1}^{n}\frac{\partial g}{\partial x_{i}}\left( X(\phi
_{1}),\ldots ,X(\phi _{n})\right) \phi _{i}.
\end{equation}

In particular, $DX(h)=h$ for every $h\in {\mathscr H}$. By iteration, one can
define the $m$-th derivative $D^{m}F$, which is an element of $L^2(\Omega ,{\mathscr H}^{\odot m})$ for every $m\geq 2$, where ${\mathscr H}^{\odot m}$ commonly stands for the $m$-th symmetric tensor product of ${\mathscr H}$. Indeed, we set

\pagebreak  

\vspace*{-28pt}  
\begin{equation}\label{iterate-Malliavin}
D^m F=\sum_{i_1,\cdots,i_m=1}^n \frac{\partial^m g}{\partial x_{i_1}\cdots \partial x_{i_m}}\Big[X(\phi_1),\cdots,X(\phi_n)\Big] \phi_{i_1}\otimes \cdots \otimes \phi_{i_m}.
\end{equation}

We stress that we rather use a symmetric tensor product instead of {the usual one} because of the celebrated Schwarz rule $\partial_x \partial y =\partial_y \partial _x$, which forces $D^m F$ to be a \emph{symmetric} element of ${\mathscr H}^{\odot m}$. For $m\geq 1$ and $p\geq 1$, ${\mathbb{D}}^{m,p}$ denotes the closure of
$\mathcal{S}$ with respect to the norm $\Vert \cdot \Vert _{m,p}$, defined by
the relation
\begin{equation*}
\Vert F\Vert_{m,p}^{p}\;=\;{\mathbb E}\left[ |F|^{p}\right] +\sum_{i=1}^{m}{\mathbb E}\left[
\Vert D^{i}F\Vert _{\mathcal{H}^{\otimes i}}^{p}\right].
\end{equation*}
We often use the notation $\mathbb{D}^{\infty} := \bigcap_{m\geq 1}
\bigcap_{p\geq 1}\mathbb{D}^{m,p}$. To justify properly the validity of the closure procedure of $\mathcal{S}$ with respect to the norm $\Vert \cdot\Vert_{m,p}$ once needs to prove that the Malliavin derivatives are \emph{closable}. Indeed, the closability is required to ensure that the limit of $DF_n$ does not depend on the choice of the approximating sequence $F_n$.
~\\\\
Another operator closely related to $D$ is the Ornstein-Uhlenbeck operator. For $F=\phi(X(h_1),\cdots,X(h_n))$, we set
\begin{equation}\label{Ornstein-uhlenbeck}
{\textbf L}[F]=\Delta \phi (X(h_1),\cdots,X(h_n))-\sum_{i=1}^n X(h_i)\frac{\partial g}{\partial x_{i}}(X(h_1),\cdots,X(h_n)).
\end{equation}

\subsection{Wiener chaos and multiple integrals}

For every $k\ge1$, we denote by $\mathcal{H}_k$ the $k$-th \textit{Wiener chaos} of $X$ defined as the closed linear subspace of $L^2(\Omega)$ generated by the family of random variables $\{H_k(X(h)), h\in {\mathscr H}, \|h\|_{{\mathscr H}}=1\}$, where $H_k$ is the $k$-th Hermite polynomial given by
\[
H_k(x)= (-1)^k e^{\frac {x^2}2} \frac {d^k} {dx^k} \left( e^{-\frac {x^2}2}\right).
\]

For any $k\geq 1$, we denote by ${\mathscr H}^{\otimes k}$ the $k$-th tensor product of ${\mathscr H}$. Set, for any $h\in{\mathscr H}$ such that $\|h\|_{\mathscr H}=1$,

\begin{equation}\label{isometry1}
I_k(h^{\otimes k})=H_k(X(h)).
\end{equation}
When $\phi=h_1^{\otimes k_1}\otimes\cdots\otimes h_p^{\otimes k_p}$ with $(h_i)_{1\le i\le p}$ an orthonormal system and $k_1+\cdots+k_p=k$, we extend (\ref{isometry1}) by
\begin{equation}\label{isometry2}
I_k(\phi)=\prod_{j=1}^p H_{k_j}(X(h_j)).
\end{equation}

{Then $I_k$ is a linear isometry between the symmetric tensor product ${\mathscr H}^{\odot k}$ (equipped with the modified norm $\sqrt{k!} \| \cdot \|_{{\mathscr H}^{\otimes k}}$) and the $k$-th Wiener chaos $\mathcal{H}_k$. In the particular case ${\mathscr H}=L^2(A,\mathcal{A}, \mu)$, where $\mu$ is a $\sigma$-finite measure without\hfilneg}

\pagebreak  

\noindent   
atoms, then $ {\mathscr H}^{\odot k}$ coincides with the space $ L^2_s(\mu^k)$  of symmetric functions which are square integrable with respect to the product measure $\mu^k$, and for any $f\in {\mathscr H}^{\odot k}$ the random variable $I_k(f)$ is commonly denoted as the multiple stochastic integral of $f$ with respect to the centered Gaussian measure generated by $X$.
~\\\\
The following fact is fundamental in the theory of Gaussian spaces.
\begin{theorem}
Any random variable $F\in L^2(\Omega)$ admits an orthogonal decomposition of the form
$$F= {\mathbb E}[F]+\sum_{k=1}^\infty I_k(f_k),$$
where the kernels $f_k\in  {\mathscr H}^{\odot k}$ are uniquely determined by $F$. In the sequel, we shall also denote $J_k(F)=I_k(f_k)$.
\end{theorem}
~\\
The random variables $I_k(f_k)$ inherit some properties from the algebraic structure of the Hermite polynomials, such that the product formula (\ref{prodfor}) below. To state it, let us introduce a definition. Let $\{e_i, i\geq 1\}$ be a complete orthonormal system in ${\mathscr H}$.
\begin{definition}\label{deficontract}
Given $f\in {\mathscr H}^{\odot k}$ and $g\in {\mathscr H}^{\odot j}$, for every $r=0,\dots, k\wedge j$, the \textit{contraction} of $f$ and $g$ of order $r$ is the element of ${\mathscr H}^{\otimes(k+j-2r)}$ defined by
\[
f\otimes_r g= \sum_{i_1, \dots, i_r=1}^\infty \langle f, e_{i_1} \otimes \cdots \otimes e_{i_r}
\rangle_{{\mathscr H}^{\otimes r}} \otimes  \langle g, e_{i_1} \otimes \cdots \otimes e_{i_r}
\rangle_{{\mathscr H}^{\otimes r}}.
\]
\end{definition}
When ${\mathscr H}=L^2({\mathbb R}_+)$, the latter formula simply becomes
\begin{align*}
 & f\otimes_r g ({\rm x},{\rm y}) \\
 & =\int_{{\mathbb R}_+^r}f(x_1,x_2,\cdots,
x_{k-r},t_1,\cdots,t_r)g(y_1,y_2,\cdots,y_{k-r},t_1,\cdots,t_r)dt_1\cdots
dt_r.
\end{align*}

The contraction $f\otimes_r g$ is not necessarily symmetric, and we denote by $f \widetilde{\otimes}_rg$ its symmetrization. We end this section by summarizing four important properties of \textit{Wiener chaos} which will be used in the sequel. {Throughout the sequel, ${\mathbf Ker}$ denotes the kernel of a linear {operator} and ${\textbf Id}$ stands for the identity operator.}
\begin{enumerate}
\item{$\mathcal{H}_k={\mathbf Ker}({\textbf L}+k{\textbf Id})$,}
\item{$\mathcal{H}_k\subset {\mathbb D}^\infty$,}
\item{For every $(p,q)\in {\mathbb N}^2$, we have the following product formula
\begin{equation}\label{prodfor}
I_p(f) I_q(g)=\sum_{r=0}^{p\wedge q} r! \binom{p}{r} \binom{q}{r} I_{p+q-2r}(f\widetilde{\otimes}_r g).
\end{equation}}
\item{For any $f\in{\mathscr H}^{\odot k}$, we have the following isometry:
$${\mathbb E}[I_k(f)^2]=k! \Vert f \Vert _{{\mathscr H}^{\otimes k}}^2.$$}
\end{enumerate}

\subsection{Main properties of Malliavin operators}

The Malliavin derivative $D$, defined in (\ref{defi-malliavin}),
obeys the following \textsl{chain rule}.
If $\varphi :\mathbb{R}^{n}\rightarrow \mathbb{R}$ is continuously
differentiable with bounded partial derivatives and if $F=(F_{1},\ldots
,F_{n})$ is a vector of elements of ${\mathbb{D}}^{1,2}$, then $\varphi
(F)\in {\mathbb{D}}^{1,2}$ and
\begin{equation}\label{e:chainrule}
D\,\varphi (F)=\sum_{i=1}^{n}\frac{\partial \varphi }{\partial x_{i}}(F)DF_{i}.
\end{equation}

The domain ${\mathbb D}^{1,2}$ can be precisely related to the Wiener chaos decomposition. Indeed,
\begin{equation}\label{Domain-D}
{\mathbb D}^{1,2}=\Big\{F\,\in L^2(\Omega)\,:\,\,\sum_{k=1}^{\infty }k\|J_k F\|^2_{L^2(\Omega)}<\infty\Big\}.
\end{equation}
In the particular case where  ${\mathscr H}=
L^{2}(A,\mathcal{A},\mu )$ (with $\mu $ non-atomic), then the
derivative of a random variable $F$ in $L^2(\Omega)$ whose chaotic expansion is
$$
F={\mathbb E}[F]+\sum_{k=1}^\infty I_k(f_k),\,\text{with}\,f_k\in{\mathscr H}^{\odot^k},
$$
can be identified with the element of $L^2(A \times \Omega )$ given by
\begin{equation}
D_{x}F=\sum_{k=1}^{\infty }kI_{k-1}\left( f_{k}(\cdot ,x)\right) ,\quad x \in A.  \label{dtf}
\end{equation}

The notation $I_{k-1}(f_k(\cdot,x))$ means that we freeze one coordinate and take the Wiener-It\^{o} integral of order $k-1$ with respect to the $k-1$ remaining coordinates. One should notice that, since the $f_k$ are taken symmetric, then the above notation do not depend on the choice of the frozen coordinate. As a matter of fact, $D_x F$ ($x\in A$) is an element of $L^2(A \times \Omega )$. We denote by $\delta $ the adjoint of the operator $D$, also called the
\textit{divergence operator}. We recall that $\delta$ exists since the operator $D$ is closed. A random element $u\in L^2(\Omega ,\mathcal{H})$
belongs to the domain of $\delta $, noted $\mathrm{Dom}\,\delta $, if and
only if it verifies
$|{\mathbb E}\left[\langle DF,u\rangle _{\mathcal{H}}|\right]\leq c_{u}\,\Vert F\Vert _{L^2(\Omega)}$
for any $F\in \mathbb{D}^{1,2}$, where $c_{u}$ is a constant depending only
on $u$. If $u\in \mathrm{Dom}\,\delta $, then the random variable $\delta (u)$
is defined by the duality relationship (customarily called \textit{integration by parts
formula})
\begin{equation}
{\mathbb E}[F\delta (u)]={\mathbb E}[\langle DF,u\rangle _{\mathcal{H}}],  \label{ipp0}
\end{equation}
which holds for every $F\in {\mathbb{D}}^{1,2}$.
More generally, if $F\in \mathbb{D}^{1,2}$ and $u\in {\rm Dom}\,\delta$ are such that the three expectations
${\mathbb E}\big[F^2\|u\|^2_{\mathscr H}]$, ${\mathbb E}\big[F^2\delta(u)^2\big]$ and
${\mathbb E}\big[\langle DF,u\rangle_{\mathscr H}^2\big]$ are finite, then
 $Fu\in{\rm Dom}\,\delta$ and
\begin{equation}\label{etoile}
\delta(Fu)=F\delta(u)- \langle DF,u\rangle_{\mathscr H}.
\end{equation}

The \textit{Ornstein-Uhlenbeck} operator ${\textbf L}$, defined in (\ref{Ornstein-uhlenbeck}) satisfies the following relation
\begin{equation}\label{DiagoL}
{\textbf L}=\sum_{k=0}^{\infty }-kJ_{k},
\end{equation}
expressing that ${\textbf L}$ is diagonalizable with spectrum $-{\mathbb N}$ with the Wiener chaos being its eigenspaces. Besides, the domain of ${\textbf L}$ is
\begin{equation}\label{domaineL}
\mathrm{Dom}({\textbf L})=\{F\in L^2(\Omega ):\sum_{k=1}^{\infty }k^{2}\left\|
J_{k}F\right\| _{L^2(\Omega )}^{2}<\infty \}=\mathbb{D}^{2,2}\text{.}
\end{equation}
There is an important relation between the operators $D$, $\delta $ and {\textbf L}.
A random variable $F$ belongs to
$\mathbb{D}^{2,2}$ if and only if $F\in \mathrm{Dom}\left( \delta D\right) $
(i.e. $F\in {\mathbb{D}}^{1,2}$ and $DF\in \mathrm{Dom}\, \delta $) and, in
this case,
\begin{equation}
\delta DF=-{\textbf L} F.  \label{k1}
\end{equation}
In particular, if $F\in\mathbb{D}^{2,2}$ and $H,G\in\mathbb{D}^{1,2}$ are such that $HG\in\mathbb{D}^{1,2}$, then
\begin{eqnarray*}
-{\mathbb E}[HG\,{\textbf L} F]&=&{\mathbb E}[HG\,\delta DF]\\
&=&{\mathbb E}[\langle D(HG),DF\rangle_{\mathscr H}]\\
&=&{\mathbb E}[H\langle DG,DF\rangle_{\mathscr H}]+{\mathbb E}[G\langle DH,DF\rangle_{\mathscr H}].\\
\end{eqnarray*}
We end this section by introducing the operator ${\textbf L}^{-1}$which is central in the next section in concern with applications of the Malliavin calculus to Stein's method. For any $F\in L^2(\Omega)$, we know that $F$ can be expanded over the Wiener chaos, namely one has $F={\mathbb E}[F]+\sum_{k=1}^\infty J_k F.$

Thus we set
$$
{\textbf L}^{-1} F={\textbf L}^{-1} \Big[ F-{\mathbb E}[F]\Big]
=\sum_{k=1}^\infty \frac {-1} {k} J_k F.
$$
One can show by using (\ref{domaineL}) that ${\textbf L}^{-1} F \in \text{Dom}({\textbf L})$.
Thus, relying on (\ref{DiagoL}) we have ${\textbf L} {\textbf L}^{-1} F =F-{\mathbb E}[F]$.
The operator ${\textbf L}^{-1}$ is called the {\it pseudo-inverse} of ${\textbf L}$.
Notice that, since $\text{Dom}({\textbf L})\subset {\mathbb D}^{1,2}$
\Big[compare (\ref{Domain-D}) and (\ref{domaineL})\Big],
for any $F\in {\mathbb D}^{1,2}$ the quantity
$<DF,-D{\textbf L}^{-1} F>_{\mathscr H}$ is well defined.

\section{Connecting Stein's method with Malliavin calculus}

As is {discussed} in Section 2, the Stein operator $L$ for normal approximation is given by $Lf(w)= f'(w) - wf(w)$ and the equation
\begin{equation}\label{Stein-id-normal}
{\mathbb E}\big[f'(Z) - Zf(Z)\big] = 0
\end{equation}
holds for absolutely continuous functions $f$ for which the expectations exist if and only if $Z \sim \mathcal{N}(0,1)$. Interestingly, this equation is nothing but a simple consequence of integration by parts. Since there is the integration by parts formula of the Malliavin calculus for functionals of general Gaussian processes, there is a natural connection between Stein's method and the Malliavin calculus. Indeed, integration by parts has been used in less general situations to construct the equation
\begin{equation}\label{Stein-general2}
{\mathbb E} [Wf(W)] = E[T f'(W)]
\end{equation}
which is Case 2 of (\ref{Stein-general}) discussed in Subsection 2.3. Let us provide two examples below.~\\

\begin{example} 
Assume ${\mathbb E} [W] = 0$ and ${\rm Var}(W) = 1$. Then we have ${\mathbb E} [T] = 1$.
If $W$ has a density $\rho > 0$ with respect to the Lebesgue measure, then by integration by parts, $W$ satisfies (\ref{Stein-general2}) with $T = h(W)$, where

$$
h(x) = \frac{\int_x^\infty y\rho(y)dy}{\rho(x)}.
$$
If $\rho$ is the density of $\mathcal{N}(0,1)$, then $h(w)= 1$ and (\ref{Stein-general2}) reduces to (\ref{Stein-id-normal}).
\end{example}
\begin{example} 
Let $X = (X_1, \dots, X_d)$ be a vector of independent Gaussian random
variables and let $g: {\mathbb R}^d \rightarrow  {\mathbb R}$ be an absolutely continuous function. Let $W = g(X)$. Chatterjee in \cite{C09} used Gaussian interpolation and integration by parts to show that $W$ satisfies (\ref{Stein-general2}) with $T = h(X)$ where
$$
h(x) = \int_0^1 \frac{1}{2\sqrt{t}}{\mathbb E} \big[\sum_{i=1}^d \frac{\partial g}{\partial x_i}(x)\frac{\partial g}{\partial x_i}(\sqrt t x + \sqrt{1 - t}X)\big]dt.
$$
If $d = 1$ and $g$ the identity function, then $W \sim \mathcal{N}(0,1)$,  $h(x)=1$, and again (\ref{Stein-general2}) reduces to (\ref{Stein-id-normal}).
\end{example}

As the previous example  shows (see Chatterjee \cite{C09} for details), it is possible to construct the function $h$ when one deals with sufficiently smooth functionals of a Gaussian vector. This is part of a general phenomenon discovered by Nourdin and Peccati in \cite{NP09}. Indeed, consider a functional $F$ of an isonormal Gaussian process $X=\{X(h), h\in {\mathscr H}\}$ over a real Hilbert space ${\mathscr H}$. Assume $F \in {\mathbb D}^{1,2}$, ${\mathbb E} [F] = 0$ and ${\rm Var}(F) = 1$.  Let $f: {\mathbb R} \rightarrow {\mathbb R}$ be  a bounded $\mathcal{C}^1$ function having a bounded derivative. Since ${\textbf L}^{-1} F \in \text{Dom}({\textbf L})$, ${\textbf L}^{-1} F \in \mathbb{D}^{2,2}$ and $D{\textbf L}^{-1}F \in \mathrm{Dom}\,\delta$.  By (\ref{k1}) and  ${\mathbb E} [F] = 0$, we have
$$
F = {\textbf L}{\textbf L}^{-1}F = \delta(-D{\textbf L}^{-1}F).
$$
Therefore
$$
{\mathbb E} [Ff(F)]={\mathbb E}[{\textbf L}{\textbf L}^{-1}F \times f(F)] = {\mathbb E}[\delta(-D{\textbf L}^{-1}F)f(F)].
$$
By the integration by parts formula (\ref{ipp0}),
$$
{\mathbb E}[\delta(-D{\textbf L}^{-1}F)f(F)] = {\mathbb E}[\langle Df(F),-D{\textbf L}^{-1}F\rangle_{\mathscr H}]
$$
and by the chain rule,
$$
{\mathbb E}[\langle Df(F),-D{\textbf L}^{-1}F\rangle_{\mathscr H}] = {\mathbb E}[f'(F)\langle DF,-D{\textbf L}^{-1}F\rangle_{\mathscr H}].
$$
Hence
$$
{\mathbb E} [F f(F)] = {\mathbb E} [f'(F)\langle DF,-D{\textbf L}^{-1}F\rangle_{\mathscr H}]
$$
and $F$ satisfies (\ref{Stein-general2}) with $T = \langle DF,-D{\textbf L}^{-1}F\rangle_{\mathscr H}$.

Now let $f_h$ be the unique bounded solution of the Stein equation (\ref{Stein-eqn-normal}) where $h: {\mathbb R} \rightarrow {\mathbb R}$ is continuous and $|h| \le 1$. Then $f_h \in \mathcal{C}^1$ and $\|f'_h\|_{\infty} \le 4\|h\|_{\infty} \le 4$,
and we have
\begin{eqnarray*}
{\mathbb E} [h(F)] - {\mathbb E} [h(Z)] &=& {\mathbb E} \{f_h'(F)[1 - \langle DF,-D{\textbf L}^{-1}F\rangle_{\mathscr H}]\} \\
&=& {\mathbb E} \{f_h'(F)[1 - {\mathbb E} (\langle DF,-D{\textbf L}^{-1}F\rangle_{\mathscr H} | F)]\}.
\end{eqnarray*}
Therefore
\begin{eqnarray*}
\sup_{h \in \mathcal{C}, |h| \le 1} |{\mathbb E} [h(F)] - {\mathbb E} [h(Z)]| &\le& \|f'_h\|_{\infty}{\mathbb E} \Big[|1 - {\mathbb E} (\langle DF,-D{\textbf L}^{-1}F\rangle_{\mathscr H} | F)|\Big] \\
&\le& 4{\mathbb E} \Big[|1 - {\mathbb E} (\langle DF,-D{\textbf L}^{-1}F\rangle_{\mathscr H} | F)|\Big].
\end{eqnarray*}
If $F$ has a density with respect to the {Lebesgue} measure, then
\begin{eqnarray*}
d_{\mathrm{TV}}({\mathscr L}(F), \mathcal{N}(0,1)) &:=& \frac{1}{2}\sup_{|h|\le 1}|{\mathbb E} [h(F)] - {\mathbb E} [h(Z)]| \\
&=& \frac{1}{2}\sup_{h \in \mathcal{C}, |h| \le 1} |{\mathbb E} [h(F)] - {\mathbb E} [h(Z)]|\\
&\le& 2{\mathbb E} \Big[|1 - {\mathbb E} (\langle DF,-D{\textbf L}^{-1}F\rangle_{\mathscr H} | F)|\Big].
\end{eqnarray*}
If, in addition, $F \in \mathbb{D}^{1,4}$, then $\langle DF,-D{\textbf L}^{-1}F\rangle_{\mathscr H}$ is square-integrable and
\begin{equation*}
{\mathbb E} \Big[|1 - {\mathbb E} (\langle DF,-D{\textbf L}^{-1}F\rangle_{\mathscr H} | F)|\Big] \le \sqrt{{\rm Var}[{\mathbb E} (\langle DF,-D{\textbf L}^{-1}F\rangle_{\mathscr H} | F)]}.
\end{equation*}
Thus we have the following theorem of Nourdin and Peccati (2011).

\begin{theorem}\label{bound-TV}
Let $F \in \mathbb{D}^{1,2}$ such that ${\mathbb E} [F] = 0$ and ${\rm Var}(F)=1$. If $F$ has a density with respect to the {Lebesgue} measure, then
\begin{equation}
d_{\mathrm{TV}}({\mathscr L}(F), \mathcal{N}(0,1)) \le 2{\mathbb E} \Big[|1 - {\mathbb E} (\langle DF,-D{\textbf L}^{-1}F \rangle_{\mathscr H} | F)|\Big].
\end{equation}
If, in addition, $F \in \mathbb{D}^{1,4}$, then
\begin{equation}\label{boundtvwiener}
d_{\mathrm{TV}}({\mathscr L}(F), \mathcal{N}(0,1)) \le 2\sqrt{{\rm Var}[{\mathbb E} (\langle DF,-D{\textbf L}^{-1}F\rangle_{\mathscr H} | F)]}.
\end{equation}
\end{theorem}

The bound (\ref{boundtvwiener}) is optimal for normal approximation for functionals of Gaussian processes. Many examples can be found in the literature, and the reader is referred to this website
\begin{verbatim}
https://sites.google.com/site/malliavinstein/home
\end{verbatim}
for a complete overview. One can also consult the good survey Peccati \cite{P14} with an emphasis on recent developments. For the sake of completeness, we shall illustrate the optimality of the bound (\ref{boundtvwiener}) using the example of the quadratic variation of a fractional Brownian motion. To do this, let $H\in(0,1)$ be the Hurst parameter of a fractional Brownian motion $\{B_t^H\}_{t>0}$. It is well known from the Breuer-Major Theorem that if $0<H<\frac{3}{4}$, then for some suitable $\sigma_H>0$,
\begin{equation*}
F_{n,H}:=\frac{1}{\sigma_H \sqrt{n}}\sum_{k=1}^n \left((B_{k+1}^H-B_k^H)^2-1\right)\xrightarrow[n\to\infty]{{\mathscr L}}~\mathcal{N}(0,1).
\end{equation*}
Similarly, if $H=\frac 3 4$, then one can prove that for some $\sigma_{\frac 3 4}>0$, we have

\begin{equation*}
F_{n,\frac 3 4}:=\frac{1}{\sigma_{\frac 3 4}\sqrt{n\log n}}\sum_{k=1}^n \left((B_{k+1}^{\frac 3 4}-B_k^{\frac 3 4})^2-1\right)\xrightarrow[n\to\infty]{{\mathscr L}}~\mathcal{N}(0,1).
\end{equation*}

Using the fact that $(B_{k+1}^H-B_k^H)^2-1$ is an element of the second Wiener chaos, and the equation (\ref{defi-malliavin}), one can obtain explicit bounds on ${\rm Var}[{\mathbb E} (\langle DF_{n,H},-D{\textbf L}^{-1}F_{n,H}\rangle_{\mathscr H}]$ (see Nourdin and Peccati \cite{NP13}).
Since
$${\rm Var}[{\mathbb E} (\langle DF,-D{\textbf L}^{-1}F\rangle_{\mathscr H} | F)]
\le {\rm Var}[{\mathbb E} (\langle DF,-D{\textbf L}^{-1}F\rangle_{\mathscr H}],$$
it follows from (\ref{boundtvwiener}) that
$$
d_{TV}(F_{n,H},G)\le c_H
\left\{\begin{array}{lrl}
\frac{1}{\sqrt{n}}&\text{if} & H\in(0,\frac 5 8)\\
\frac{(\log n)^{\frac 3 2}}{\sqrt{n}}&\text{if}& H=\frac 5 8\\
n^{4H-3}&\text{if}& H\in (\frac 5 8,\frac 3 4)\\
\frac {1}{\log n}&\text{if}& H=\frac 3 4.
\end{array}
\right.
$$

These bounds are shown to be optimal in Nourdin and Peccati \cite{NP13}. We wish to mention that for $H>\frac 3 4$, $F_{n,H}$ does not converge to a Gaussian distribution. Instead, it converges to the so-called Hermite distribution, which belongs to the second Wiener chaos and is therefore not Gaussian.

%
%
%

\section{The Nualart-Peccati criterion of the fourth moment and Ledoux's idea}

\subsection{Some history}

During the year 2005, in the seminal article \cite{NP05}, David Nualart and Giovanni Peccati discovered the following {remarkable} fact. Take $F_n=I_k(f_n)$ a sequence of elements of the $k$-th Wiener chaos. Then, $F_n$ converges in law towards the Gaussian measure $\mathcal{N}(0,1)$ \textit{if and only if} ${\mathbb E}[F_n^2]\to 1$ and ${\mathbb E}[F_n^4]\to 3$. This result can be seen as a drastic simplification of the so-called \textit{method of moments} which consists in proving that {${\mathbb E}[F^{p}]\to\int_{\mathbb R} x^{p}e^{-\frac{x^2}{2}}\frac{dx}{\sqrt{2\pi}}$},\, \textit{for each positive integer p.} A bit later, in  Peccati and Tudor \cite{PT04}, the strength of this theorem was considerably reinforced by its multivariate counterpart. Indeed, Giovanni Peccati and Ciprian Tudor proved that a random vector $F_n=(F_{1,n},\cdots,F_{d,n})$ with chaotic components converges in distribution towards a Gaussian vector with covariance $C$ \textit{if and only if} the covariance of $F_n$ converges to $C$ and for each $i\in\{1,\cdots,d\}$, $F_{i,n}$ converges in distribution to $\mathcal{N}\big[0,C(i,i)\big]$. That is to say, for chaotic random variables, the \textit{componentwise} converges implies the \textit{joint} convergence. This observation yielded to a very efficient strategy for proving central convergence in the Wiener space, by somehow \textit{decomposing} the convergence on each Wiener chaos. This approach, which is more and more used by practitioners as an alternative to the semi-martingale approach, is at the heart of a very active line of research. See this website for details and an exhaustive overview.
\begin{verbatim}https://sites.google.com/site/malliavinstein/home\end{verbatim}

\subsection{Overview of the proof of Nourdin and Peccati}

Whereas the original proof of the fourth moment Theorem relied on some tools of stochastic analysis, Ivan Nourdin and Giovanni Peccati produced a new proof in \cite{NP09} based on a suitable combination of Malliavin calculus and Stein's method. As noticed by the Nourdin and Peccati in \cite{NP09} and in conclusion of previous sections 3 and 4,
one is left to show that
\begin{equation*}
{\rm Var}({\mathbb E}[T_1|W])\le {\rm Var}\Big[\langle DF_n,-D {\textbf L}^{-1} F_n\rangle_{\mathscr H}\big]\to 0,
\end{equation*}
under the assumptions that ${\mathbb E}[F_n^4]\to 3$ and ${\mathbb E}[F_n^2]\to 1$, if $\{F_n\}$ is a sequence of elements in the $k$-th Wiener chaos. In fact, following their strategy we will prove that
\begin{eqnarray}\label{Fundament-equation}
{\rm Var}\Big[\langle DF_n,-D {\textbf L}^{-1} F_n\rangle_{\mathscr H}\Big]\le \frac{k-1}{3 k} \Big[{\mathbb E}[F_n^4]-3{\mathbb E}[F_n^2]\Big].
\end{eqnarray}

\underline{Step 1: computing ${\rm Var}\big[\langle DF_n,-D {\textbf L}^{-1} F_n\rangle_{\mathscr H}\big]$.}~\\

By equation (\ref{dtf}), one has $DF_n(t)=DI_k(f_n)(t)=k I_{k-1}(f_n(\cdot,t))$. By applying the product formula for multiple integrals (\ref{prodfor}), we get:
\begin{eqnarray*}
\frac{1}{k}\langle D F_n, D F_n\rangle_{\mathscr H}&=& k \int_0^\infty I_{k-1}(f_n(\cdot,t))^2 dt\\
&=&k \int_0^\infty \sum_{r=0}^{k-1}r! \binom{k-1}{r}^2 I_{2k-2-2r}\big[f_n(\cdot,t)\widetilde{\otimes}_r f_n(\cdot,t)\big]dt\\
&=&k\sum_{r=0}^{k-1}r! \binom{k-1}{r}^2 I_{2k-2-2r}\Big[\int_0^\infty f_n(\cdot,t)\widetilde{\otimes}_r f_n(\cdot,t)dt\Big]\\
&=&k \sum_{r=1}^{k} (r-1)! \binom{k-1}{r-1}^2I_{2k-2r}\Big[f_n\widetilde{\otimes}_r f_n\Big].\\
&=&k \sum_{r=1}^{k-1} (r-1)! \binom{k-1}{r-1}^2I_{2k-2r}\Big[f_n\widetilde{\otimes}_r f_n\Big]+ k!\Vert f\Vert_{{\mathscr H}}^2.
\end{eqnarray*}

Taking into account that ${\mathbb E}[\langle D F_n, D F_n\rangle_{\mathscr H}]=k {\mathbb E}[F_n^2]=k k!\Vert f\Vert_{{\mathscr H}}^2,$ the orthogonality of the Wiener chaos entails that:
\begin{equation}{\label{Step1}}
{\rm Var}\big[\langle DF_n,-D {\textbf L}^{-1} F_n\rangle_{\mathscr H}\big]=\sum_{r=1}^{k-1}\frac{r^2}{k^2}(r!)^2\binom{k}{r}^4(2k-2r)!\Vert f_n\widetilde{\otimes}_r f_n\Vert_{{\mathscr H}^{\odot (2k-2r)}}^2.
\end{equation}

\underline{Step 2: computing ${\mathbb E}[F_n^4]$.}~\\

By product formula again, we have

$$F_n^2=\sum_{r=0}^k r! \binom{k}{r}^2 I_{2k-2r}(f_n\widetilde{\otimes}_r f_n).$$
This yields to
$${\mathbb E}[F_n^4]=\sum_{r=0}^k (r!)^2 \binom{k}{r}^4 (2k-2r)! \Vert f_n\widetilde{\otimes}_r f_n\Vert_{{\mathscr H}^{\odot (2k-2r)}}^2.$$
Unfortunately the latter expression is not immediately comparable with (\ref{Step1}) because the ``zero contractions" $\Vert f_n\widetilde{\otimes}_0 f_n\Vert_{{\mathscr H}^{\odot (2k-2r)}}^2$ do not appear in the decomposition of $\langle D F_n, D F_n\rangle_{\mathscr H}$. To face this problem, one has to take into account an intermediary expression, namely ${\mathbb E}[F_n^2 \langle D F_n,  D F_n\rangle_{\mathscr H}]=\frac{k}{3} {\mathbb E}[F_n^4]$. After few combinatorial arguments respective to the symmetrizations of the contractions of the kernels $f_n$ appearing in the scalar product, and by a suitable comparison of the terms ${\mathbb E}[F_n^4]$, ${\rm Var}(\langle D F_n , D F_n\rangle_{\mathscr H})$ and ${\mathbb E}[F_n^2 \langle D F_n , D F_n\rangle_{\mathscr H}]$ one may show that:
\begin{equation}\label{Step2}
{\mathbb E}[F_n^4]-3{\mathbb E}[F_n^2]^2=\frac{3}{k}\sum_{r=1}^{k-1}r (r!)^2 \binom{k}{r}^4 (2k-2r)!\Vert f_n \widetilde{\otimes}_{2k-2r} f_n \Vert_{{\mathscr H}^{\odot (2k-2r)}}^2.
\end{equation}
Since this part is rather technical and irrelevant for our purpose, we refer the reader to the book Nourdin Peccati \cite{N12} for more precise arguments. By comparing the equations (\ref{Step1}) and (\ref{Step2}), one may recover (\ref{Fundament-equation}).

\subsection{About Ledoux's generalization}

As we just showed, the usual way of proving the inequality (\ref{Fundament-equation}) relies heavily on the product formula for multiple Wiener integrals. In particular, this gives the impression that the various combinatorics coefficients appearing in the formulae (\ref{Fundament-equation}), (\ref{Step1}) and (\ref{Step2}) are playing a major role in the phenomenon. In 2012, in the very insightful article \cite{le}, Michel Ledoux tried to tackle this problem by using a more ``algebraic" approach based on spectral theory and $\Gamma$-calculus. In particular, he was able to prove the inequality (\ref{Fundament-equation}) without using the product formula, under some suitable spectral conditions. Unfortunately, the provided conditions seemed rather difficult to check in practice, and no new structure (in addition to the Wiener space one) with a fourth moment phenomenon was given.

\section{The  general fourth moment Theorem for Dirichlet forms}\label{general}

In this section, we keep the remarkable intuition of Ledoux in \cite{le} of exploiting the algebraic and spectral properties of the chaotic random variables rather than product formulae techniques. However, at the very beginning, we will take a different path. As we will show, the fourth moment Theorem turns out to be a direct consequence of a very simple spectral assumption. This simplification will allow us to produce various examples of new structures where the phenomenon holds. We will make a crucial use of the very powerful formalism of Dirichlet forms to achieve our goals.

\subsection{The Dirichlet structures}

Originated from potential theory and physics, the Dirichlet form theory has become a central object in analysis. Under the impulsion of Beurling, Deny, Fukushima, Meyer or Mokobodzki, just to name a few, the Dirichlet forms theory has unveiled beautiful ramifications in many area of mathematics such as geometric measure theory, partial differential equations, Markov processes, Malliavin calculus... Here, we will adopt a modern formalism, namely the notion of Dirichlet structure.~\\

We will say that a Dirichlet structure is a set $(E,\mathcal{F},\mu,{\mathbb D},\Gamma)$ such that:
\begin{enumerate}
\item{$(E,\mathcal{F},\mu)$ is a probability space,}
\item{${\mathbb D}$ is a vector space \textit{dense} in $L^2(\mu)$,}
\item{$\Gamma: {\mathbb D} \times {\mathbb D} \rightarrow L^1(\mu)$ is a \textit{bilinear, symmetric} and \textit{non-negative} operator,}
\item{${\mathbb D}$ endowed with the norm $\Vert X \Vert_{{\mathbb D}}=\sqrt{{\mathbb E}\big[X^2\big]+{\mathbb E}\big[\Gamma[X,X]\big]}$ is \textit{complete}.}
\item{For any vectors $X=(X_1,\cdots,X_p)\in{\mathbb D}^p$ and $Y=(Y_1,\cdots,Y_q)\in{\mathbb D}^q$, for any functions $F\in\mathcal{C}^1({\mathbb R}^p,{\mathbb R})\cap\text{Lip}$ and $G\in\mathcal{C}^1({\mathbb R}^q,{\mathbb R})\cap\text{Lip}$  we have $F(X)\in {\mathbb D}$ and $G(Y)\in{\mathbb D}$. Besides, one has the following~~\textit{functional calculus}
\begin{equation}\label{functional calculus}
\Gamma[F(X),G(Y)]=\sum_{i=1}^p\sum_{j=1}^q \partial_i F(X) \partial_j G(Y) \Gamma[X_i,Y_j].
\end{equation}}
\end{enumerate}

The mapping $(X,Y)\rightarrow \mathcal{E}\big[X,Y\big]={\mathbb E}\big[\Gamma[X,Y]\big]$ is customarily called the \textit{Dirichlet form} associated {with} the \textit{carr\'{e}-du-champ} operator $\Gamma$ . For any of theses structures, one may associate another operator ${\textbf L}$ defined on some domain $\text{dom}({\textbf L})$ \text{dense} in ${\mathbb D}$ such that the next integration by parts holds for any $(X,Y)\in \text{dom}({\textbf L}) \times \text{dom}({\textbf L})$.
\begin{equation}\label{IPPDiri}
{\mathbb E}\big[\Gamma[X,Y]\big]=-{\mathbb E}\big[X {\textbf L}[Y]\big]=-{\mathbb E}\big[Y{\textbf L}[X]\big].
\end{equation}
Besides, one may derive from (1-5) the next relation between $\Gamma$ and ${\textbf L}$.
\begin{equation}\label{GL}
2\Gamma[X,Y]={\textbf L}[XY]-Y{\textbf L}[X]-X{\textbf L}[Y]
\end{equation}
In addition to assumptions (1-5) which characterize a general Dirichlet structure, we assume further that
\begin{itemlist}[H00]
\item[(H1)]{$-{\textbf L}$ is \textit{diagonalizable} with \textit{spectrum} $\{0=\lambda_0<\lambda_1<\lambda_2<\cdots<\lambda_p<\cdots\}$,}
\item[(H2)] For each $X\in{\mathbf Ker}({\textbf L}+\lambda_k {\textbf Id})$, $$X^2\in\bigoplus_{\alpha\le 2\lambda_k}{\mathbf Ker}({\textbf L}+\alpha {\textbf Id}).$$
\end{itemlist}
Before stating and proving our result, we stress that the previous framework covers the case of the Wiener structure. Indeed, for the Wiener structure, $\Gamma[F,G]=\langle DF, DG\rangle_{\mathscr H}$ defined on ${\mathbb D}^{1,2}$ and the corresponding ${\textbf L}$ is the Ornstein-Uhlenbeck operator. Properties (1-3) are straightforward, (4) proceeds from the closability of ${\mathbb D}^{1,2}$ see subsection (3.2), and (5) is a consequence of the chain rule (\ref{e:chainrule}). (H1) comes from (\ref{DiagoL}) and (H2) is a trivial consequence of the product formula (\ref{prodfor}).

\subsection{Fourth moment theorem for Dirichlet structures with (H1) and (H2)}

Let $X$ be an eigenfunction of $-{\textbf L}$ associated with eigenvalue $\lambda_k$ with ${\mathbb E}[X^2]=1$. We will show that
\begin{equation}\label{Goal}
{\rm Var}(\Gamma[X,X])\le \frac{\lambda_k^2}{3}  ({\mathbb E}[X^4]-3)
\end{equation}
Suppose that the inequality (\ref{Goal}) is true and take $\phi$ any test function. On the one hand we have by the chain rule (5) that $\Gamma[\phi(X),X]=\phi'(X)\Gamma[X,X]$. On the other hand, by integrations by parts (\ref{IPPDiri}) we also have ${\mathbb E}[\phi'(X)\Gamma[X,X]]=\lambda_k{\mathbb E}[X\phi(X)]$. Thus we are in the setting of the case 2 of the section 2.3 and one has the inequality
\begin{equation}\label{basic}
d_{\mathrm{TV}}({\mathscr L}(X),\mathcal{N}(0,1))\le {\mathbb E}\Big|\frac{\Gamma[X,X]}{\lambda_k}-1\Big|\le\frac{1}{\lambda_k}\sqrt{{\rm Var}(\Gamma[X,X]).}
\end{equation}
We used above the fact that ${\mathbb E}\big[\Gamma[X,X]\big]=-{\mathbb E}\big[X{\textbf L}[X]\big]=\lambda_k$. Now, using (\ref{GL}) and taking into account that $X$ is actually an \textit{eigenfuntion} of ${\textbf L}$, we derive $\Gamma[X,X]=\frac{1}{2}\big[{\textbf L}+2\lambda_k{\textbf Id}\big][X^2]$. Hence, we get $\big[\Gamma[X,X]-\lambda_k\big]=\frac{1}{2}\big[{\textbf L}+2\lambda_k{\textbf Id}\big][X^2-1]$. The rest of the proof is as follows. We use (H2) to say that $X^2$ and hence $X^2-1$ are expanded in finitely many eigenspaces of $-{\textbf L}$ with maximum eigenvalue being $2\lambda_k$.
\begin{eqnarray*}
{\rm Var}(\Gamma[X,X])&=&{\mathbb E}\big[\frac{1}{2}\big[{\textbf L}+2\lambda_k{\textbf Id}\big][X^2-1]\times \frac{1}{2}\big[{\textbf L}+2\lambda_k{\textbf Id}\big][X^2-1]\big]\\
&=&\frac{1}{2}{\mathbb E}\big[{\textbf L}[X^2-1]\times \frac{1}{2}\big[{\textbf L}+2\lambda_k{\textbf Id}\big][X^2-1]\big]\\
&+&\lambda_k {\mathbb E}\big[(X^2-1)\times \frac{1}{2}\big[{\textbf L}+2\lambda_k{\textbf Id}\big][X^2-1]\big]\\
&\underset{(H2)}{=}&\frac{1}{4}\sum_{\lambda_i\le 2\lambda_k}(-\lambda_i)(2\lambda_k-\lambda_i){\mathbb E}\big[J_k[X^2-1]^2\big]\\
&+& \lambda_k{\mathbb E}\Big[(X^2-1) \big[\Gamma[X,X]-\lambda_k\big]\Big]\\
&\underset{(H2)}{\le}& \lambda_k {\mathbb E}\Big[(X^2-1)\Gamma[X,X]\Big]-\lambda_k^2 {\mathbb E}[X^2-1]\\
&\underset{(5)}{=}& \lambda_k{\mathbb E}\Big[\Gamma[\frac{X^3}{3}-X,X]\Big]\\
&\underset{(\ref{IPPDiri})}{=}&\lambda_k^2\Big[\frac{{\mathbb E}[X^4]}{3}-1\Big].\\
\end{eqnarray*}
So, the proof of (\ref{Goal}) is done and one has shown (by dividing by $\lambda_k$) that

\begin{equation}\label{VarGamma_1}
\text{Var}(\Gamma[X,\frac{X}{\lambda_k}])\leq \Big[\frac{{\mathbb E}[X^4]}{3}-1\Big].
\end{equation}
Relying on the equations (\ref{functional calculus}) and (\ref{IPPDiri}) it is easy to see that for any function $f$ bounded with bounded derivatives

$${\mathbb E}\left[f'(X) \Gamma[X,\frac{X}{\lambda_k}]\right]={\mathbb E}\left[f(X) X\right].$$

Since $\Gamma[X,\frac{X}{\lambda_k}]$ plays the same role as the term ``$\langle DF,-DL^{-1} F\rangle_{{\mathscr H}}$" of Section 4, Theorem \ref{bound-TV} applies and can be combined with (\ref{VarGamma_1}). Now take $X_n$ as a sequence of eigenfunctions of ${\textbf L}$ with same eigenvalue, such that ${\mathbb E}[X_n^2]=1$. We deduce that
$$d_{\mathrm{TV}}({\mathscr L}(X_n),\mathcal{N}(0,1))\le \sqrt{\frac{{\mathbb E}[X_n^4]}{3}-1}.$$
In particular, if ${\mathbb E}[X_n^4]\to 3$, then we have convergence to $\mathcal{N}(0,1)$ in total variation.

\subsection{Dirichlet structures with (H1) and (H2)}

Here we give two examples of Dirichlet structures which satisfy (H1) and (H2). The first and most important example is provided by the Wiener space. First from equation \ref{DiagoL}, one can see that the so-called Ornstein-Uhlenbeck operator is diagonalizable with spectrum ${\mathbb N}$ and the Wiener chaos are its eigenspaces. This fact shows that the assumption (H1) is fulfilled. One the other hand, take $k\in{\mathbb N}$ and $X=I_k(f)$ some element of the $k$-th Wiener chaos. By using the product formula (\ref{prodfor}),
$$X^2=\sum_{r=0}^k k! \binom{k}{r}^2 I_{2k-2r}(f\tilde{\otimes}_r f).$$

From which one can deduce that $$X^2\in\bigoplus_{r\le 2k}{\mathbf Ker}({\textbf L}+r{\textbf Id}).$$
Hence (H1) and (H2) hold true for the Wiener structure. One should notice, that the precise combinatorial coefficients arising in the product formula play no role in the assumption (H2). In some sense, the approach by Dirichlet forms being less combinatoric is more direct than the usual proof in Nualart Peccati \cite{NP05}.~\\

Not only the Dirichlet forms approach simplifies the proof of the Nualart-Peccati criterion, but also enables one to give new examples of Dirichlet structures where a fourth moment phenomenon holds. In order to avoid technicalities, we restrict the exposition to the finite dimensional case. So, let us introduce the Laguerre Dirichlet structure. Let $\nu \geq -1$, and $ \pi_{1,\nu}(dx) =
x^{\nu-1}\frac{e^{-x}}{\Gamma(\nu)} \textbf{1}_{(0,\infty)}(x) d x$ be the Gamma
distribution with parameter $\nu$ on ${\mathbb R}_+$. The associated Laguerre generator is
defined for any test function $\phi$ (in dimension one) by:
\begin{equation}\label{Lag1} {\textbf L}_{1,\nu} \phi(x)= x\phi''(x)+(\nu+1-x)\phi'(x).
\end{equation}
By a classical tensorization procedure, we obtain the Laguerre generator in
dimension $d$ associated to the measure $ \pi_{d,\nu}(d x) = \pi_{1,\nu}(d x_1)
\pi_{1,\nu}(d x_2) \cdots \pi_{1,\nu}(d x_d)$, where $x=(x_1,x_2, \cdots,x_d)$.
\begin{equation}\label{Lag2} {\textbf L}_{d,\nu}\phi(x) = \sum_{i=1}^{d}
\Big{(}x_{i} \partial_{i,i}\phi+(\nu+1-x_i)\partial_i \phi\Big{)}
\end{equation}

It is well known that (see for example~Ledoux \cite{le}) that the spectrum of
${\textbf L}_{d,\nu}$ is given by~$-{\mathbb N}$ and moreover that
\begin{equation} \textbf{Ker}({\textbf L}_{d,\nu} + p{\textbf Id}) =
\left\{\sum_{i_1+i_2+\cdots+i_{d}=p} \alpha(i_1,\cdots,i_{d})\prod_{j=1}^{d}
L^{(\nu)}_{i_j}(x_j)\right\},
\end{equation}
where $L^{(\nu)}_n$ stands for the Laguerre polynomial of order $n$ with parameter $\nu$ which is defined by
$$ L_n^{(\nu)}(x)= \frac{x^{-\nu} e^x}{n!} \frac{d^n}{dx^n}
\left(e^{-x} x^{n+\nu}\right).$$
Again, we have the following decomposition:
\begin{equation}\label{decompo2} L^2({\mathbb R}^d,\pi_{d,\nu})=\bigoplus_{p=0}^\infty
\textbf{Ker}({\textbf L}_{d,\nu} + p {\textbf Id})
\end{equation}

As a matter of fact, assumption (H1) {follows}. Let us check assumption (H2). Assume now that $X$ is an eigenfunction of
${\textbf L}_{d,\nu}$ with eigenvalue $-\lambda_p=-p$. In particular, $X$ is a
multivariate polynomial of degree $p$. Therefore,  $X^2$ is a multivariate
polynomial of degree $2p$. Note that by expanding $X^2$ over the basis of
multivariate Laguerre polynomials $\prod_{j=1}^{d} L^{(\nu)}_{i_j}(x_j), i_{j}\ge
0$, we get that $X^2$ has a finite expansion over the first eigenspaces of the
generator ${\textbf L}_{d,\nu}$, i.e.
\begin{equation*} X^2 \in \bigoplus_{p=0}^{M} {\mathbf Ker}({\textbf L}_{d,\nu}+ p {\textbf Id}).
\end{equation*} Finally, by taking care of degree reasons we infer that $M=2p$ and thus the assumption (H2) is valid in this structure.

\pagebreak

\enlargethispage*{15pt}  
\section{Acknowledgments}

We would like to thank Giovanni Peccati for his valuable comments which have helped to improve the exposition of this paper. This work is partially supported by AFR grant 4897114 at the University of Luxembourg and by Grant C-146-000-034-001 and Grant R-146-000-182-112 at the National University of Singapore.

\end{document}